\documentclass[11 pt]{amsart}
\usepackage{amssymb}
\usepackage{amsmath}
\usepackage{amsfonts}
\usepackage{graphicx}
\usepackage{amsthm}
\usepackage{enumerate}
\usepackage[mathscr]{eucal}
\usepackage{verbatim}

 \theoremstyle{plain}
\newtheorem{theorem}{Theorem}
\newtheorem{lemma}[theorem]{Lemma}

\theoremstyle{plain}

\theoremstyle{remark}

\newcommand{\Conv}{\mathop{\scalebox{2.0}{\raisebox{-0.2ex}{$\ast$}}}}

\def\bbR{{\mathbb {R}}}

\begin{document}

\date{April, 2012}

\title{Some toy Furstenberg sets and projections of the four-corner Cantor set}

\author[]{Daniel M. Oberlin}

\address
{D. M.  Oberlin \\
Department of Mathematics \\ Florida State University \\
 Tallahassee, FL 32306}
\email{oberlin@math.fsu.edu}

\subjclass{28A75}
\keywords{Hausdorff dimension, Furstenberg set}


\begin{abstract}
We give lower bounds for the Hausdorff dimensions of some model Furstenberg sets. 
\end{abstract}

\maketitle

In \cite{W} Wolff noted that the following question stems from work of Furstenberg: fix
$\alpha\in (0,1)$ and consider the class of compact sets $E \subset \bbR^2$ 
which have the property that for each direction $\theta$ in $\bbR^2$ there is a line in the direction of $\theta$ which intersects $E$ in a set of Hausdorff dimension at least $\alpha$. What is the minimum possible Hausdorff dimension $\dim (E )$ for such an $E$? Wolff showed in \cite{W} that this minimum must lie in the interval 

\begin{equation}\label{a}
[\max \{1/2+\alpha  , 2\alpha \},1/2+3\alpha /2].
\end{equation}
 The only subsequent progress
concerns the case $\alpha =1/2$ (where the lower bounds $1/2+\alpha$ and $2\alpha$ coincide): 
in this case there is some $\epsilon >0$ such $\dim (E)\geq 1+\epsilon $. 
This is a result of Bourgain and Katz-Tao. Specifically, it is a consequence of
Bourgain's work \cite{B} on the Erd\H os-Volkmann ring conjecture (about the existence of subrings of $\bbR$ having 
Hausdorff dimension strictly between $0$ and $1$) combined with work of Katz and Tao \cite{KT} on the equivalence
of special cases of the ring conjecture, of the Furstenberg problem, and of Falconer's distance problem.

The main purpose of this note is to give an improvement of the lower bound of \eqref{a} for a narrow class of specific, but rather natural, examples of these Furstenberg sets $E$. To describe our $E$'s we begin by recalling Kahane's 
construction \cite{K} of a Besicovitch set (see also \cite{SS}): let $C$ be the $1/2$-dimensional Cantor set
$$
\big\{\sum_{j=1}^\infty \epsilon_j 4^{-j}:\epsilon_j =0,3\big\}.
$$
For an ordered pair $(c_1 ,c_2 )\in C\times C$, the so-called four-corner Cantor set, let $\ell_{c_1 c_2}$ be the line segment in $\bbR^2$ joining the points 
$(0,c_1)$ and $(1,c_2 /2)$. Let $B$ be the union of all such segments $\ell_{c_1 c_2}$. Then $B$ is a compact
subset of $\bbR^2$. Since the slope of $\ell_{c_1 c_2}$ is $c_2 /2-c_1$ and since 
$$C/2-C=C/2+C-1=[-1,1/2],
$$ 
it follows that 
$B$ contains a line segment with slope $m$ for each $m\in [-1,1/2]$. For $0\le x\le 1$, the $x$-section of $B$ is the 
set 
\begin{equation}\label{b}
\{(1-x)c_1 +xc_2 /2 :(c_1 ,c_2 )\in C\times C\}, 
\end{equation}
a projection into $\bbR$ of the $1$-dimensional and purely unrectifiable Cantor 
set $C\times C$. It follows from the Besicovitch projection theorem (or see \cite{SS} for an elementary proof) that,
for almost all $x\in\bbR$, the set \eqref{b} has $1$-dimensional Lebesgue measure $0$. Therefore $B$ has $2$-dimensional Lebesgue measure $0$. Thus the union of some finite collection of rotations of $B$ is a Besicovitch set - a compact subset of the plane having $2$-dimensional Lebesgue measure $0$ and containing a unit line segment in each direction. Now suppose that $0<\alpha <1$ and that $K\subset [0,1]$ is any compact set satisfying 
$\dim(K)=\alpha$. 
(In this note the term {\it dimension} and the function $\dim$ will always refer to Hausdorff dimension.)
Then, for $(c_1 ,c_2 )\in C\times C$, the set $\{(1-x)c_1 +xc_2 /2 :x\in K\}\subset\ell_{c_1 ,c_2}$ has dimension 
$\alpha$ and so, if  
\begin{equation}\label{c}
E_\alpha =\{(1-x)c_1 +xc_2 /2 :(c_1 ,c_2 )\in C\times C, x\in K\}, 
\end{equation}
it follows that for each $m\in [-1,1/2]$ there is some line segment 
$\ell_{c_1 c_2}$
with slope $m$ which intersects $E_\alpha$ 
in a set of Hausdorff dimension at least $\alpha$. The set $E_\alpha$ is similar to the set which Wolff constructs in \cite{W} to give the upper bound in \eqref{a}. 
In fact, Wolff's example can easily be modified to yield sets $E_\alpha$ 
as in \eqref{c} with $\dim (E_\alpha )\le 1/2+3\alpha /2$. 
We would like to prove that $\dim (E_\alpha )\ge 1/2+3\alpha /2$, i.e., that the so-called Furstenberg Conjecture is true 
at least for the sets $E_\alpha$. But
we have only the following lower bound for $\dim (E_\alpha )$:
\begin{theorem}\label{result}
If $E_\alpha$ is as in \eqref{c} (with $dim(K)\ge \alpha$), then 
\begin{equation}\label{d}
dim (E_\alpha )\ge c+3\alpha /2 
\end{equation}
where $c=\log (8/3)/\log (16)\approx .35$. 
%
%
\end{theorem}
\noindent 
Note that 
\eqref{d} improves the lower bound in \eqref{a} (but only for our particular $E_\alpha$'s)
whenever $1/2-c<\alpha < 2c$. 

The next result follows immediately from Theorem 5.8 in \cite{F}. It provides a natural approach 
to \eqref{d}.
\begin{lemma}\label{lemma}
Suppose that $0<s,t\le 1$ and that $E\subset \bbR^2$ is such that the $x$-sections of $E$ 
have Hausdorff dimension at least $t$ for an $x$-set having Hausdorff dimension at least $s$. Then 
$\dim (E)\ge s+t$.
\end{lemma}
\noindent Thus we write $E_{\alpha ,x}$ for the $x$-section of $E_{\alpha}$ and note that
\begin{equation}\label{dim}
\dim(K)\ge \alpha \text{ and } \dim (E_{\alpha ,x})\ge \tau \ \forall x\in K 
\Rightarrow  \dim (E_\alpha )\ge \alpha +\tau .
\end{equation}
Now, for $x\in K$, \eqref{b} shows that $E_{\alpha ,x}$ is the projection $P_x (C\times C)$ where 
%
$P_x :\bbR^2 \rightarrow \bbR$
%
is given by $P_x (x_1 ,x_2 )=(1-x)x_1 +xx_2 /2$.
(Projections of $C\times C$ were studied from another point of view in \cite{NPV} and in \cite{BV} - see also \cite{BLV}.)
Our aim here is to apply estimates on the dimensions of projections of $C\times C$ in order to obtain results like \eqref{d}.
As an easy example of this strategy, note that 
$$
\dim (E_{\alpha ,x})=\dim \big( P_x (C\times C)\big) \ge 1/2
$$
for any $x\in K$ (since 
$\dim (C)=1/2$). With $\dim (K)=\alpha$, \eqref{dim} gives the $1/2+\alpha$ lower bound from \eqref{a}
for $\dim (E_\alpha )$.
Again, here is a theorem due to Kaufman \cite{Kf} which refines an earlier result of Marstrand: 
\begin{theorem}\label{K}
Suppose the compact set
$A\subset\bbR^2$ has dimension $\beta \leq 1$ and $p_\omega (A)$ is the projection of $A$ onto the line through the origin in the direction of $\omega\in S^1$, then 
\begin{equation}\label{dimest}
\dim\{\omega\in S^1 :\dim p_\omega (A) <\gamma\}\leq \gamma 
\end{equation}
whenever $\gamma\leq\beta$.
\end{theorem}
\noindent With the correspondences
$$
\omega =\frac{(1-x,x/2)}{\sqrt{(1-x)^2 +(x/2)^2}},\ p_\omega (A) =\frac{1}{\sqrt{(1-x)^2 +(x/2)^2}}\, P_x (A)
$$
and the fact that $\dim (C\times C)=1$, Kaufman's result 
and $\dim (K)=\alpha$
show that,  for any $\epsilon >0$, $\dim (E_{\alpha ,x})\ge \alpha -\epsilon$ must hold for $x$ in an 
$\alpha$-dimensional subset of $K$. With \eqref{dim} this recovers,
for $E_\alpha$, the $2\alpha$ lower bound of \eqref{a}. To go further we will improve the conclusion of 
Theorem \ref{K} when $A=C\times C$. The paper \cite{O} contains the conjecture that 
the conclusion \eqref{dimest} of Kaufman's theorem is, in fact, always improvable to  
\begin{equation}\label{conjecture}
\dim\{\omega\in S^1 :\dim p_\omega (A) <(\beta +\gamma)/2\}\leq \gamma . 
\end{equation}
For $\gamma =\beta$ this follows from \eqref{dimest}, while for $\gamma =0$ this is proved in \cite{O}. 
(In \cite{B2} Bourgain proves that given $\gamma\in (0,1)$ there is $\kappa =\kappa (\gamma )>\beta /2$
such that 
\begin{equation*}
\dim\{\omega\in S^1 :\dim p_\omega (A) <\kappa\}\leq \gamma . )
\end{equation*}
Applied with $A=C\times C$, $\beta =1$, and $\gamma =\alpha -\epsilon$, the conjectural conclusion \eqref{conjecture} would imply that 
\begin{equation}\label{f}
\dim\{x\in \bbR :\dim P_x (C\times C) <(1+\alpha -\epsilon)/2\}\leq \alpha -\epsilon
\end{equation}
for small $\epsilon$. 
(The modification of Wolff's example mentioned after \eqref{c} shows that \eqref{f}, if true, would be sharp,
and it follows from Lemma \ref{lemma2} below that we do at least have
\begin{equation*}
\dim\{x\in\bbR :\dim P_x (C\times C) <c+\alpha /2 -\epsilon \}\leq \alpha -\epsilon .)
\end{equation*}
Now, since $P_x (C\times C)=E_{\alpha ,x}$
and $\dim (K)=\alpha$, it would follow from \eqref{f} 
that $\dim (E_{\alpha ,x}\ge (1+\alpha -\epsilon )/2$ for an $\alpha$-dimensional set 
of $x$'s in $K$ and so, by \eqref{dim}, 
that $\dim (E_{\alpha })\ge 1/2+3\alpha /2$, the best possible result.
But even when $A$ is a product set (like $C\times C$)
the conjecture \eqref{conjecture} is likely to be difficult: 
such a conclusion would imply 
the aforementioned Erd\H os-Volkmann ring conjecture in the case of
rings $R$ of dimension $\sigma$ not exceeding $1/2$.
To see this we again reindex projections and, for $t\in\bbR$, write $P^t (x_1 ,x_2 )=x_1 +tx_2$. If we assume that 
$\dim (R)=\sigma$ and take $A=R\times R$, the conclusion \eqref{conjecture} would yield (with $\beta =2\sigma , \,
\gamma =2\epsilon >0$)
\begin{equation*}
\dim\{t\in \bbR :\dim (R+tR) <\sigma +\epsilon\}\leq 2\epsilon .
\end{equation*}
This is impossible if $R+R\cdot R\subset R$ and $2\epsilon <\sigma$.

We return to Theorem \ref{result}.
Suppose $0<\alpha ' <\alpha$ and let $\mu$ be a nonnegative compactly-supported 
Borel measure on $K$ satisfying $\mu (I)\lesssim |I|^{\alpha '}$ for 
all intervals $I\subset \bbR$. We will show that for any $\epsilon >0$ we have
\begin{equation*}
\dim \big( P_x (C\times C) \big)\ge c+\alpha ' /2-\epsilon
\end{equation*}
for a set of $x$'s having full $\mu$-measure. Since such an $x$-set must have dimension at least $\alpha '$,
\eqref{d} will then follow from Lemma \ref{lemma}. Since it is more convenient to work 
with the projections $P^t$ instead of $P_x$,
we will actually prove the following lemma.

%
%
%
%
\begin{lemma}\label{lemma2}
Suppose that $0<\alpha ' <1$ and that $\mu$ is a compactly-supported 
probability measure on $\bbR$ satisfying $\mu (I)\lesssim |I|^{\alpha '}$ for 
all intervals $I\subset \bbR$. Then for any $\epsilon >0$ we have
\begin{equation}\label{g}
\dim ( C+tC )\ge c+\alpha ' /2-\epsilon
\end{equation}
for a Borel set of $t$'s having full $\mu$-measure.
\end{lemma}

\noindent To establish \eqref{g}, let $\lambda$ be the Cantor-Lebesgue measure 
%
\[
\Conv\nolimits_{k=1}^\infty  \Big(\frac{1}{2}\big(\delta_0 +\delta_{3\cdot 4^{-k}}\big)\Big)
\]
%
on $C$ and, for $t\in\bbR$, write $\lambda_t$ for the dilate of $\lambda$ supported on $tC$.
Then $\lambda\ast\lambda_t$ is supported on $C+tC$ and so \eqref{g} will follow from 
\begin{equation*}
\int_\bbR \int_\bbR \frac{|\hat{\lambda}(s)\hat{\lambda}(ts)|^2}{|s|^{1-\tau}}\,ds\, d\mu (t)<\infty
\end{equation*}
whenever $\tau <c+\alpha ' /2$. This will be a consequence of the estimate
\begin{equation}\label{mainest2}
\int_{4^K}^{4^{K+1}}|\hat{\lambda}(s)|^2 \,\int_\bbR |\hat{\lambda}(ts)|^2 \, d\mu (t)\, ds\lesssim \big(
4^{1-c-\alpha ' /2}\big)^K .
\end{equation}
Now 
$$
 \hat{\lambda}(s)=\prod_{k=1}^\infty \Big(\frac{1}{2}\big(1+e^{-2\pi i 3\cdot 4^{-k}s}\big)\Big)
 ,\ \  |\hat{\lambda}(s)|^2 
=\prod_{k=1}^\infty \cos^2 (3\pi 4^{-k}s) .
$$
With 
$$
P(s)=\prod_{k=1}^\infty \cos^2 ( 2\pi  4^{-k}s) ,
$$
the bound \eqref{mainest2} is a consequence of the two estimates 
\begin{equation}\label{Pest1}
\int_{4^K}^{4^{K+1}}P(s)\, ds\lesssim 4^{K/2}
\end{equation}
and 
\begin{equation}\label{Pest2}
\int_\bbR P(ts)\, d\mu (t)\lesssim  \big(6^{1/2}4^{-(1+\alpha ')/2}\big)^K \text{ if } 
4^K\le s\le 4^{K+2}
\end{equation}
along with the fact that $c=\log (8/3)/\log (16)$ so that $4^{1-c}=6^{1/2}$.
The proofs of  \eqref{Pest1} and \eqref{Pest2} are not difficult. To see \eqref{Pest1}  we write
\begin{multline*}
\int_{4^K}^{4^{K+1}}P(s)\, ds=4^K \int_1^4 \prod_{k=1}^\infty \cos^2 ( 2\pi  4^{K-k}s)\, ds\leq \\
4^K \int_1^4 \prod_{k=1}^K \cos^2 ( 2\pi  4^{K-k}s)\, ds =3\cdot 4^K \int_0^1 P_K (s)\, ds
\end{multline*}
where 
\begin{equation}\label{k}
P_K (s)= \prod_{l=0}^{K-1} \cos^2 ( 2\pi 4^{l}s) =
\frac{1}{4^{K}}\big|\prod_{l=0}^{K-1} \big(e^{-2\pi i4^l s}+e^{2\pi i4^l s}\big)
\big|^2 
\end{equation}
so that
\begin{equation*}
 \int_0^1 P_K (s)\, ds=\frac{2^{K}}{4^{K}}.
\end{equation*}

To establish \eqref{Pest2} we begin by
%
%
%
%
fixing an even Schwartz function $\rho$ with 
\begin{equation}\label{rho}
|\rho |\lesssim \sum_{j=0}^\infty 2^{-2j}\chi_{[-2^{j},2^j ]},\ \hat{\rho}=1 \text{ on }[-100,100]. 
\end{equation}
Write $s=4^K s_0$ with $1\le  s_0 \le 16$ and
define the $\alpha '$-dimensional measure $\tilde{\mu}$ by 
$$
\int f\, d\tilde{\mu}=\int f(ts_0 )\, d\mu (t).
$$
Then \eqref{Pest2} will follow as above from 
\begin{equation}\label{est2}
\int P_K (t)\, d\tilde{\mu}(t) \lesssim \big(6^{1/2}4^{-(1+\alpha ' )/2}\big)^K .
\end{equation}
Expanding the product in \eqref{k} shows that 
$$
P_K (t)= \frac{1}{4^{K}} 
\sum_{\epsilon_k ,\epsilon_k ' =\pm 1}e^{2\pi i \sum_{k=0}^{K-1} (\epsilon_k +\epsilon_k ')4^k t}.
$$
Let $\rho_K$ be the dilate of $\rho$ defined by $\rho_K (t)=4^K \rho (4^K t)$. Now $P_K$ is a trigonometric polynomial and, by \eqref{rho}, $\widehat{\rho_K}=1$ on the support of the discrete measure
$\widehat{P_K}$. Thus $P_K=P_K \ast \rho_K$ and so 
\begin{equation}\label{est3}
\int P_K (t)\, d\tilde{\mu}(t)=\int P_K (t) \, \big(\rho_K \ast \tilde{\mu}(t)\big)\, dt .
\end{equation}
To bound \eqref{est3} we will estimate certain $L^2$ norms of $P_K$ and $\rho_K \ast \tilde{\mu}$. To begin we
compute $\| P_K \|_{L^2 ([0,1]}$.
Let $S_K$ be the set of integers $n$ which can be represented
\begin{equation}\label{S_K}
n=\sum_{k=0}^{K-1} \delta_k (n)\,4^k ,\ \delta_k (n) \in \{-2,0,2\}
\end{equation} 
and for $n\in S_K$ let $r(n)$ be the cardinality of the set
\begin{equation*}
\Big\{ (\epsilon_0 ,\dots ,\epsilon_{K-1} ,\epsilon_0 ' \dots ,\epsilon_{K-1} ' ):
\epsilon_k ,\epsilon_k ' \in \{-1,1\},\
n= \sum_{k=0}^{K-1} (\epsilon_k +\epsilon_k ')4^k \Big\}.
\end{equation*}
Then
\begin{equation*}\label{P_K}
\| P_K \|^2_{L^2 ([0,1]}=\frac{1}{4^{2K}}\sum_{n\in S_K}r(n)^2 .
\end{equation*}
Since any $n\in S_K$ has a unique representation \eqref{S_K}, $r(n)=2^{s(n)}$ where $s(n)$ is the number of indices 
$k$ in $n$'s representation \eqref{S_K} for which $\delta_k (n)=0$. Also, for any $I\subset \{0,\dots K-1 \}$ 
of cardinality $s$, there are $2^{K-s}$ elements $n\in S_K$ for which $\delta_k (n)=0$ precisely when $k\in I$. Thus 
\begin{equation}\label{P_K2}
\| P_K \|^2_{L^2 ([0,1]}=\frac{1}{4^{2K}}\sum_{s=0}^K \binom {K}{s} 2^{K-s}\, 2^{2s} =\frac{2^{K}\,3^{K}}{4^{2K}}.
\end{equation}
To estimate $\rho_K \ast\tilde{\mu}$ we note that 
$\|\rho_K \ast\tilde{\mu}\|_{L^\infty (\bbR )}\lesssim 4^{K(1-\alpha ')}$ follows from \eqref{rho} and the 
estimate $\tilde{\mu}(I)\lesssim |I|^{\alpha '}$ for intervals $I$. Thus the trivial estimate 
$\|\rho_K \ast\tilde{\mu}\|_{L^1 (\bbR )}\lesssim 1$ implies 
$\|\rho_K \ast\tilde{\mu}\|_{L^2 (\bbR )}\lesssim 4^{K(1-\alpha ' )/2}$. If $\rho_K \ast\tilde{\mu}$ were compactly supported then this $L^2$ estimate together with \eqref{est3} and \eqref{P_K2} and the periodicity of $P_K$ would imply \eqref{est2}. To deal with the tail which arises because $\rho_K \ast\tilde{\mu}$ is not compactly supported, choose $M$ so that $\tilde{\mu}$ is supported in $[-M/2,M/2]$.
We will estimate 
$$
\int_{\{|t|>M\}} P_K (t) \, \big(\rho_K \ast \tilde{\mu}(t)\big)\, dt
$$
by writing $\{|t|>M\}$ as an essentially disjoint  union of intervals $[n,n+1]$ and noting that the estimate
$|\rho_K \ast\tilde{\mu}|\lesssim 4^{-K} n^{-2}$ on such an $[n,n+1]$ 
follows from $$|\rho_K (x)\lesssim 4^{-K}x^2$$ (which itself is a consequence of \eqref{rho}). Thus 
we can estimate
$$
\Big|\int_{[n,n+1]} P_k (t) \,\big(\rho_K \ast\tilde{\mu}(t)\big)\, dt \Big|\lesssim \frac{\| P_K \|_{L^2 ([0,1]}}{4^K n^2}.
$$
Summing these estimates on $n$ and taking account of \eqref{P_K2} gives 
$$
\int_{\{|t|>M\}} P_K (t) \, \big(\rho_K \ast \tilde{\mu}(t)\big)\, dt
\lesssim \big(6^{1/2}4^{-2}\big)^K
$$
and so \eqref{est2} is established.


We close with two observations: the reader may have noted that if one could improve \eqref{est2} to
\begin{equation}\label{eq10}
\int P_K (t)\, d\tilde{\mu}(t) \lesssim \big(2^{-\alpha}\big)^K
\end{equation}
for an $\alpha$-dimensional measure $\tilde{\mu}$, then \eqref{f} would follow. But one can use self-similarity to show 
that 
\begin{equation*}
\int P_K (t)\, d\lambda (t)\gtrsim (3/4)^K
\end{equation*}
so that \eqref{eq10} must fail at least for $\alpha =1/2$. On the other hand, 
there are various ways to obtain marginal improvements of our lower bound for $\dim (E_\alpha )$.
For example, 
it is possible to estimate higher order $L^p$ norms of $P_K$ and to use these estimates,
as above, to obtain further improvements of the lower bounds for $\dim (E_\alpha )$. In particular, 
\begin{equation*}
\| P_K \|_{L^3 ([0,1])}=\Big(\frac{5}{16}\Big)^{(K+1)/3}.
\end{equation*}
This leads to 
\begin{equation*}
\dim (E_\alpha )\ge c' +4\alpha /3 ,\  c'=\frac{\log (32/5)}{6\log 2}
\end{equation*}
which improves \eqref{d} for some values of $\alpha <2c$ and also improves \eqref{a} for certain values 
of $\alpha < 1/2-c$. But such incremental improvements do not seem worth pursuing.


\end{document}